# Speed Trajectory Planning at Signalized Intersections Using Sequential Convex Optimization

Xianan Huang, Huei Peng

*Abstract*— An algorithm is developed to optimize vehicle speed trajectory over multiple signalized intersections with known traffic signal information to minimize fuel consumption and travel time, and to meet ride comfort requirements using sequential convex optimization method. A comparison between the proposed method and dynamic programming is carried out to verify its optimality. In addition, vehicle motion during turning is studied because of its significant effect on fuel consumption and travel time.

## I. INTRODUCTION

In city driving, stop-and-go and idling due to congestion and signalized intersections cause significantly increased fuel consumption. A Recent study [1] shows that on average fuel consumption at signalized intersections accounts for more than 50% of the total consumption for a whole trip. The main techniques used to address this issue are adaptive traffic signal controls such as SCOOT [2] and SCATS [3]. These infrastructure centric solutions have limitations, however, due to the responsiveness of traffic flow to traffic signals and the reduced effectiveness when the number of vehicles is low.

Dedicated short-range communication (DSRC) supports short range and reliable data communication between vehicles and infrastructures [4, 5], which enables vehicle centric approach. Together with connected automated vehicle (CAV) technologies, the vehicle centric solution can improve fuel efficiency, mobility and safety. Broadcast by road side equipment (RSE) through DSRC, signal phase and timing (SPaT) contains the current and future signal phase and timing information, enabling predictive control and smooth driving. The NHTSA performed a preliminary analysis on the benefits of SPaT, showing a 90% reduction in red light violations and up to 35% of savings in energy [6].

Information from traffic signals available via SPaT allows the vehicle speed trajectory to be planned so as to reduce fuel consumption at signalized intersections, a concept known as eco-approach/departure. Multi-stage optimization methods have been used to solve the optimization problem [7, 8]. With the goal of avoiding stops at signalized intersections, vehicle speed is controlled at the maximum speed without having to stop at intersections. Xia et al. [8] experimentally studied the effect of speed advisory with rule-based speed planning, and found a 14% reduction in fuel consumption and a 1% reduction in travel time. Subsequently, with the smoothened speed profile designed from simplified rules, a variety of optimal trajectory following methods are used. Asadi et al. [7] has used a model predictive control algorithm with the objective function defined as a weighted sum of the trajectory following error and fuel consumption. The work has been extended to vehicle platoons [9] and hybrid electric vehicles [10]. However, the planned speed trajectory is based primarily on avoiding idling at intersections. The potential in fuel savings is not fully addressed.

To realize the full potential in fuel savings, optimization methods such as dynamic programming technique [11] and Pontryagin's minimum principle [12] were used to optimize speed trajectory in the whole problem horizon with single intersection, no other vehicles, and precisely known traffic signal assumptions. Inaccurate traffic signal states in the problem horizon has been extended by [13]. [14] has included other queuing vehicles at the intersection in the analysis. They estimated the queue clearing time and used pseudospectral method to obtain the optimal speed trajectory. In [15], a discretized solution using Dijkstra's algorithm has been obtained for multiple intersections. They assumed the vehicles cross the intersection only at a specific point in time such as at the beginning, middle or end of a green phase window. The discrete choices were modeled as nodes in a graph and solved a shortest path problem. In many of the works cited above, additional assumptions, for example, constant traveling speed along each road section, are made to reduce computation load. In addition, though, because urban driving entails frequent turns, vehicles may incur significant a penalty in terms of both fuel economy and time. To the best of our knowledge, however, the effect of turning has not been considered in the literature.

In this paper, we present a speed trajectory optimization algorithm with turning motion constraints using the sequential convex optimization method [16]. Sequential convex optimization is a method for obtaining a local optimal solution by forming convex sub-problems sequentially. It finds a local optimal solution in a timely manner without suffering from the curse of dimensionality. Another benefit is its flexibility to the form of objective function due to the sequential convex approximation. We assume that the traffic signal state is known within the problem horizon, and we do not consider the influence of a lead vehicle. The problem can then be solved over the whole problem horizon, thus taking advantage of the full potential of speed variation over multiple sections of the road. The second advantage is that we do consider turning at intersections. The turning speed constraint is determined by considering the characteristics of the intersection. The third advantage lies in the flexibility of the proposed method: it can consider multiple objectives,

The authors are with the Department of Mechanical Engineering, University of Michigan, Ann Arbor, Michigan 48109-2133, USA (e-mail: xnhuang@umich.edu, hpeng@umich.edu).

and can be applied to multiple-vehicle and multiple-intersection cases. In addition to having a flexible problem formulation, it is also important to use a robust numerical solver. We use Gurobi as the solver [17] here.

The rest of the paper is organized as follows. The model of a passenger car is constructed in Section 2. Section 3 explains the optimization problem. Section 4 presents the optimization results and the analysis. Finally, conclusions are given in Section 5.

## II. VEHICLE MODEL

### A. Fuel Consumption Model

In this study, we consider a passenger car equipped with a 4-cylinder 2.5-liter internal combustion engine and a continuously variable transmission (CVT). A simplified powertrain model is used to focus on the study of eco-approach/departure, which is possible through the following assumptions: (1) powertrain efficiency is simplified to a static look-up table; (2) CVT keeps the engine operating along the best brake specific fuel consumption (BSFC) line; (3) a simple longitudinal dynamics of the vehicle [18] is used.

$$M\dot{v} = F - Mg\sin\theta - Mgf\cos\theta - 0.5\rho C_d A(v+v_w)^2 \quad (1)$$

where $M$ is vehicle mass, $v$ is vehicle speed, $F$ is the longitudinal force, $g$ is the gravity coefficient, $\theta$ is the road grade, $f$ is the rolling resistance coefficient, $\rho$ is the air density, $C_d$ is the drag coefficient, $A$ is the vehicle cross-sectional area, and $v_w$ is the wind speed. In the following, we assume flat road and zero wind speed. The driving force is a function of gear ratio and engine torque

$$F = i_g i_f \eta_T / r_w \quad (2)$$

where $i_g$ is the transmission gear ratio, $i_f$ is the final drive ratio, $\eta_T$ is the transmission efficiency, $r_w$ is the wheel radius. The fuel consumption is estimated from the static fuel consumption map, showing in Fig. 1.

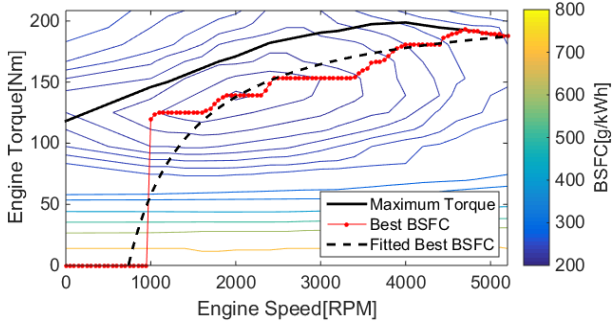

Fig. 1 BSFC map of engine model

The idling engine speed is 800 RPM and the idling torque is assumed to be 0 Nm. We assume engine stop-start technique is not available, so idling fuel consumption is not zero. The optimal BSFC point is around 2000 RPM, 140 Nm. To incorporate the transient effect of engine operation on fuel consumption, we follow the method of Li et al. [19] by adding a modification term to the fuel consumption obtained from the static fuel consumption map. The total fuel consumption is

$$Q = Q_{static} + k_e \dot{T}_e \quad (3)$$

where $Q$ is the total fuel consumption, $Q_{static}$ is the fuel consumption rate from the static lookup table, $k_e$ is the coefficient for engine transient operations, $T_e$ is the engine torque. The coefficient $k_e$ is obtained from the driving cycle FTP-72 assuming that the transient engine operations increase the fuel consumption by 4~5% [19].

The transmission is assumed to be controlled so that the engine stays on the best BSFC line

$$\omega_{opt} = \frac{b}{1-kT_{opt}} \quad (4)$$

where $\omega_{opt}$ is the engine speed along the best BSFC line, $T_{opt}$ is the engine torque, $k$ and $b$ and parameters to be identified. Under this ideal CVT assumption, the fuel consumption rate is a function of the engine power.

### B. Effect of Turning

We assume turning mainly imposes speed and acceleration limits when the vehicle is moving through an intersection. For the speed constraint, we consider the simplified unbanked surface turn model [20], which computes speed limit due to friction limit

$$v_{max} = \sqrt{Rg\mu} \quad (5)$$

where $R$ is the turning radius, $\mu$ is the friction coefficient. In addition, we assume there is a limit on vehicle speed during turning due to ride comfort

$$v = \sqrt{Ra_y} \quad (6)$$

### C. Baseline Driver Deceleration/Acceleration Model

To evaluate the effectiveness of the algorithm, a human behavior model at intersection is used as the benchmark. The deceleration and acceleration behavior model shown in Eq. (7) is from [21], which was confirmed to match experimental data very well.

$$a = ra_m\theta(1-\theta^m)^2 \quad (7)$$

In Eq. (7), $ra_m$ is function of $m$, $\theta$ is the relative acceleration/deceleration time, defined as time divided by desired acceleration/deceleration time. The model parameters are all adopted from [21]. The reaction distance is defined as the maximum distance to the intersection where driver starts to decelerate if the light state is red. The desired acceleration /deceleration time and distance is calculated with empirical functions from [21]

$$t_a = \frac{v_f - v_i}{0.5778 + 0.0669(v_f - v_i)^{1/2} - 0.0182v_i} \quad (8)$$

$$x_a = (0.467 + 0.0072v_f - 0.0076v_i)(v_f + v_i)t_a \quad (9)$$

$$t_d = \frac{x_d}{(0.473 + 0.0056v_i - 0.0049v_f)(v_f + v_i)} \quad (10)$$

$$ra_m = [2(m+1)(m+2)/m^2]|v_f - v_i|/t_{a(d)} \quad (11)$$

where $x_a$ and $x_d$ are desired acceleration/deceleration distance, $t_a$ and $t_d$ are desired acceleration/deceleration time, $v_f$ and $v_i$ are desired final and initial speed. The reaction distance is set as 150 m to the intersection. We assume that the desired deceleration distance is the distance to the intersection when the light is red and the driver is within the reaction distance. If the speed limit is 17.88m/s, $m$ for deceleration is -0.7193, and 9.1244 for acceleration from (8)(9)(10) and [21].

III. ECO-DRIVING PROBLEM AND SOLUTION METHODOLOGY

A. Mixed-Integer Problem Formulation

The speed trajectory optimization problem is formulated as a non-convex optimization problem. The objective is to minimize fuel consumption and traveling time while meet ride comfort requirement over the planning horizon; the constraints include speed limits, acceleration limit, and red light violation. The vehicle motion is discretized with sampling time, and during each sampling time, acceleration is assumed to be constant. In discrete time model, speed and displacement are

$$v(k+1) = v(k) + a\Delta t \quad (12)$$

$$d(k+1) = d(k) + \frac{(v(k+1) + v(k))}{2}\Delta t \quad (13)$$

The traction power at each time step is derived from the longitudinal vehicle model (1).

$$P(k) = Ma(k)v(k) + Mgfv(k) + 0.5\rho C l_d A v(k)^3 \quad (14)$$

As discussed in Section 2, fuel consumption is only a function of the engine power along the BSFC line. Therefore, fuel consumption $FC(k)$ is

$$FC(k) = C_f(P_{eng})P_{eng}(k) = C_f(P_{eng})P(k)/\eta_T \quad (15)$$

where $C_f(P_{eng})$ is the fuel consumption coefficient.

A travel time penalty is imposed through a negative vehicle speed term over the planning horizon, and a penalty on acceleration and jerk represents the desire for better ride comfort.

$$J_{comfort}(k) = a(k)^2 + w_j(a(k) - a(k-1))^2 \quad (16)$$

The final objective function is defined as a weighted sum of fuel consumption, traveling time and ride comfort.

$$J = w_{fc}\sum_{k=1}^{T} J_{fc}(k) - w_t \frac{1}{T}\sum_{k=1}^{T} v(t) + w_c \sum_{k=1}^{T} J_{comfort}(k) \quad (17)$$

where $T$ is the horizon time, $w_{fc}$, $w_t$, $w_c$ are weighting parameter for fuel consumption, traveling time and ride comfort respectively.

To en sure the vehicle crosses the intersection without violating the red light, we define the constraints with respect to green phase window. $t_{r2g}^{(i)}$ is defined as the time the light changes from red to green for the $i$-th green phase window of the subject intersection, and $t_{g2r}^{(i)}$ is defined as the time the light changes from green to red. These time steps are the critical times for speed trajectory optimization at signalized intersections. To put the constraints into a matrix form, we define the vehicle location at the critical times and the indicator of crossing windows as follows

$$\mathbf{k} = [k_1,...,k_N]^T \quad (18)$$

$$\mathbf{d_{r2g}} = [d_{r2g}^{(1)},...,d_{r2g}^{(N)}]^T \quad (19)$$

$$\mathbf{d_{g2r}} = [d_{g2r}^{(1)},...,d_{g2r}^{(N)}]^T \quad (20)$$

where $\mathbf{k}$ is a singleton vector with only the indicator of a selected crossing green light equals to 1, and the other elements set to 0. $N$ is the total number of green phase windows in the planning horizon at the subject intersection. $\mathbf{d}_{r2g}$ and $\mathbf{d}_{g2r}$ are vectors of vehicle locations at the critical times. With the variables defined in the vector form, the constraint for valid intersection crossing can be defined as

$$\sum_{i=1}^{N} k_i = 1, k_i \in \{0,1\} \quad (21)$$

$$\mathbf{k^T d}_{r2g} < 0 \quad (22)$$

$$\mathbf{k^T d}_{g2r} > 0 \quad (23)$$

Other constraints include the speed limit constraint $v_{max}$, the acceleration limit constraint and jerk constraint. Unlike the study performed by[8], we do not allow the vehicle to exceed the speed limit to catch a green light.

$$v(k) \in [0, v_{max}], a(k) \in [a_{min}, a_{max}] \quad (24)$$

$$a(k) - a(k-1) \in [Jerk_{min}, Jerk_{max}] \quad (25)$$

As discussed above, the problem is formulated as a non-convex optimization problem, with speed and displacement as the state variables, and acceleration and the crossing green phase window indicator as the input variables. Among the variables the crossing green phase window indicator is an integer variable. The constraints are either linear or quadratic. However, the objective function is non-convex, with nonlinear fuel consumption and aero-resistance. To solve the problem, sequential convex optimization is applied. Sequential convex optimization finds a local optimal solution by forming a convex sub-problem of the original problem sequentially. The method has been used to solve trajectory planning for aircraft, manipulator and humanoid robot [16, 22]. To make the approximation at each iteration valid, the trust region method is applied, that is, an additional constraint is applied to make the step size small. At each iteration, the two non-convex terms are approximated by the values from the previous iteration. At iteration $j+1$ the objective function of fuel consumption is shown below, which is formulated in a symmetric form.

$$J_{fc}^{j+1} = diag(fc_k^j / p_k^j)(\mathbf{a}^T(M\mathbf{D} \\ + 0.5\rho C_d A \mathbf{D}^T diag(v_k^j)\mathbf{D})\mathbf{a} + \mathbf{v}_0^T M\mathbf{a} \quad (26) \\ + \rho C_d A \mathbf{v}_0^T diag(v_k^j)\mathbf{D}\mathbf{a} + Mgf\mathbf{D}\mathbf{a} + K)$$

where $K$ is a constant term related only to the initial speed, $p_k^j$ is the traction power at iteration $j$ time $k$, $fc_k^j$ is the fuel consumption rate at iteration $j$ time $k$, $\mathbf{a}$ is the vector form of the acceleration in the planning horizon, $\mathbf{D}$ is a N×N lower triangle matrix representing the kinematic model (12)

$$\mathbf{D} = \begin{pmatrix} 1 & 0 & \cdots & 0 \\ 1 & 1 & \cdots & 0 \\ \vdots & \vdots & \ddots & \vdots \\ 1 & 1 & \cdots & 1 \end{pmatrix}, \mathbf{v} = \mathbf{Da} + v_0 \quad (27)$$

The assumption here is that in the trust region, the fuel consumption and the speed of the last iteration are valid approximations of the true value. The trust region method would impose additional linear constraint on speed and acceleration

$$v(k)^{j+1} \in [v(k)^j - \rho_v, v(k)^j + \rho_v] \quad (28)$$

$$a(k)^{j+1} \in [a(k)^j - \rho_a, a(k)^j + \rho_a] \quad (29)$$

where $\rho_v$ and $\rho_a$ are the trust region radius of speed and acceleration respectively. It is also noted from the solver that since breach-and-bound is used to solve the mixed-integer problem, the application of the trust region at each iteration would reduce the size of the search tree.

Since the multi-objective optimization problem is solved by the weighted sum method, the objective function is not guaranteed to be positive-semidefinite; thus, at each iteration standard sequential quadratic programming (SQP) is used to obtain the solution.

To initialize the sequential convex optimization, the initial cost function for fuel consumption is set to minimize the traction power rather than minimizing fuel consumption.

$$J_{fc}^0 = \mathbf{a}^T \mathbf{MDa} + \mathbf{v}_0^T \mathbf{Ma} \tag{30}$$

### B. Incorporation of Turning Motion

As discussed in the previous section, we assume the geometry of the intersection can be neglected when incorporating the turning motion, which is modeled as speed and acceleration limits as follows

$$0 \leq v(t_{cross}) \leq v_{turn} \tag{31}$$

$$a_{turn\_min} \leq a(t_{cross}) \leq a_{turn\_max} \tag{32}$$

where $t_{cross}$ is the crossing time, $v_{turn}$ is the maximum speed during turning, $a_{turn\_min}$ and $a_{turn\_max}$ are acceleration limits determined by the intersection. Since the intersection crossing time is unknown even when the crossing window is determined, the crossing speed and acceleration constraints are achieved through soft constraints, as shown in Fig. 2 and Fig. 3.

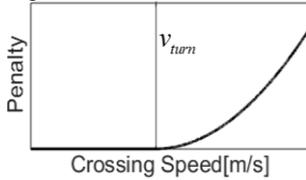
Fig. 2 Crossing Speed Soft Constraint

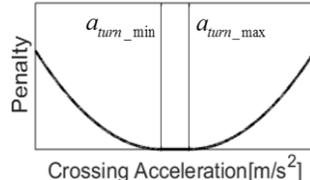
Fig. 3 Crossing Acceleration Soft Constraint

The soft constraint is implemented with a piecewise linear objective in Gurobi [17]. With the convex nature of the quadratic and the SQP approximation of the original problem, using soft constraints would preserve the convexity of the sub-problem at each iteration. However, due to the usage of the trust region, at each iteration the converging step size is small. In addition, the usage of soft constraints increases computation time for the mixed-integer programming. When the crossing time change between two consecutive iterations is larger than a specific threshold, we reinitialize the sequential convex optimization by resetting the cost function as propelling power, removing the trust region constraint and adding linear constraints for crossing speed and acceleration. We can thus achieve a closer start point to the local optimal solution.

To define the stopping criteria for the sequential optimization, the distance of improvement between the iterations is defined. The criterial iteration variables are fuel consumption rate, vehicle speed, and crossing speed. The distance of improvement is defined as

$$\Delta fc^j = \max_k \left( | fc(k)^j - fc(k)^{j-1} | \right) \tag{33}$$

$$\Delta v^j = \max_k \left( | v(k)^j - v(k)^{j-1} | \right) \tag{34}$$

$$\Delta v_{cross}^{\ j} = \left( v(t_{cross})^j - v_{turn} \right)_{v(t_{cross})^j > v_{turn}} \tag{35}$$

$$\Delta G^j = \sqrt{\Delta fc^{j2} + \Delta v^{j2} + \Delta v_{cross}^{\ j2}} \tag{36}$$

where $\Delta G^j$ is the total difference between two consecutive iterations evaluated at iteration $j$, defined as the square root of the sum of the squares of the difference in fuel consumption rate, vehicle speed and crossing speed. The iteration stops when the total difference is less than the threshold.

## IV. OPTIMIZATION RESULTS AND DISCUSSION

We first start from a single vehicle, single intersection case. The problem horizon is set to be 90 seconds. The speed limit is 17.88 m/s, or 40 mph. The acceleration limits are ±3 m/s$^2$, as used in [7]. The jerk limits are set to be ±0.5 m/s$^3$. Mixed integer programming is known to be NP hard and the solving time is related to the number of integer states and the problem size. The problem is usually solved with branch-and-bound[17]. For our case, the integer variable is the crossing window indicator and the number is small in the problem horizon. The problem is solved with a computer with Intel i7-4710MQ CPU and 16 G RAM. When the turning motion is not considered, the solving time varies from 0.4 s to 1.9 s depending on the traffic light status. When the turning motion is considered, the solving time increases dramatically, varying from 6.6 s to 8.4 s depending on traffic light status and the gap between initial speed and the desired turning speed.

### A. Checking of Optimality

Sequential convex optimization is a method for obtaining local optimal solutions of non-convex problems. To verify the optimality of the solution, the speed trajectory is compared with solutions from dynamic programming for sanity check. Although dynamic programming obtains global optimal solutions solving the problem backwards in time, the algorithm is computationally expensive and suffers from

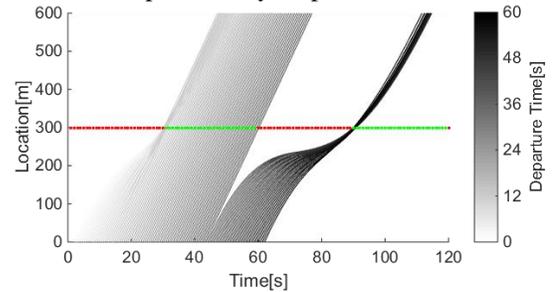
Fig. 4 Vehicle Trajectories for Different Signal Phase

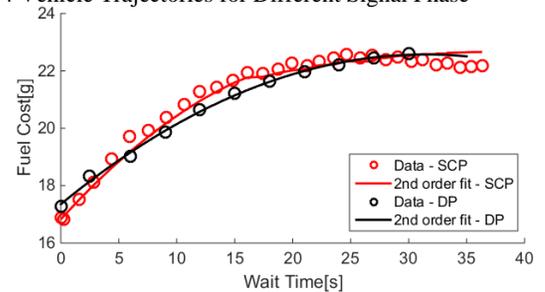
Fig. 5 Comparison between DP and SCP

curse of dimensionality. With pre-computed cost-to-go, DP would take 628 s to obtain the optimal solution. The time weight is set to be 2000 (m/s)$^{-1}$. The speed trajectories for different traffic light phase is shown in Fig. 4, with red dots representing the red phase of the traffic light and green dots representing the green phase. The change in the signal phase is achieved through fixing the traffic signal and changing the vehicle departure time.

Fig. 5 shows results from both DP and SCP. The waiting time is defined as the time difference between the actual travel time and free flow travel time; fuel cost is defined as the fuel consumed from 300 m before the intersection to 300 m after the intersection and reaching the original speed. The relation between fuel cost and wait time can be fitted with a 2$^{nd}$ order curve. The maximum difference between results from DP and SCP is 4.28%.

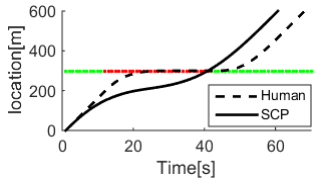
Fig. 6 Sample Comparison of Location with Driver Model

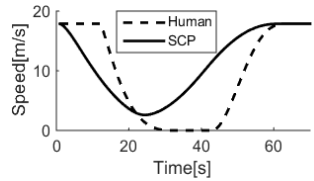
Fig. 7 Sample Comparison of Speed with Driver Model

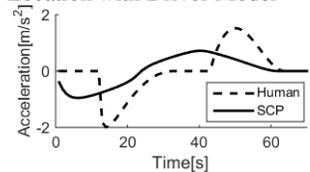
Fig. 8 Sample Comparison of Acceleration with Driver Model

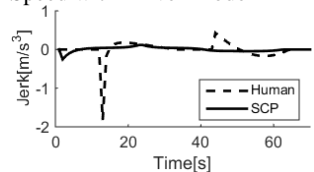
Fig. 9 Sample Comparison of Jerk with Driver Model

The optimal results are obtained for different traffic signal phases. The average fuel consumption reduction is 12.1% and time reduction is 7.5% for single intersection cases compared with the driver model. The reduction in fuel can be up to 35.6% and the reduction in time can be 16.4% depending on the traffic signal status. A sample trajectory comparison of location, speed, acceleration and jerk are shown in Fig. 6-Fig. 9. The optimization results show smoother behavior compared with the results from a driver model, mainly due to the fact that future traffic light status is known.

*B. Turning Motion Consideration*

To address the benefit of including turning motion in the optimization, we consider a left turn at an intersection with

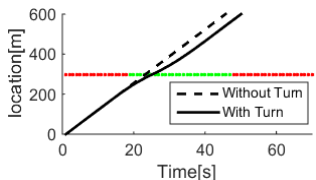
Fig. 10 Sample Location Trajectory during Turning

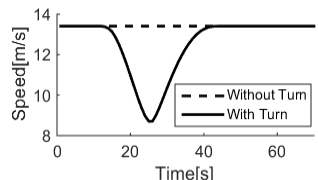
Fig. 11 Sample Speed Trajectory during Turning

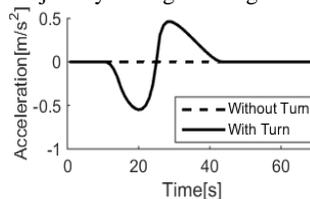
Fig. 12 Sample Acceleration Trajectory during Turning

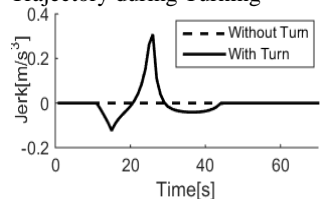
Fig. 13 Sample Jerk Trajectory during Turning

four-lane roads on each side. The turning radius is set to be 25 m, the comfort lateral acceleration level is set to be 3 m/s$^2$, and road friction coefficient is 0.7. For this study, we set the longitudinal acceleration in the turning to be 0. The maximum speed to pass through the intersection is 13.1 m/s from (5), and the comfortable maximum speed to pass through the intersection is 8.7 m/s from (6). The study is carried out with different signal phase. If the speed limit of the road is higher than the maximum safe passing speed, the method without turning motion constraints cannot obtain a feasible solution for free flow since the optimal solution is passing the intersection at constant speed. We set the speed limit to be 13 m/s for a valid comparison. A sample of a comparison between optimal trajectories with and without turning motion consideration is shown in Fig. 10-Fig. **13**.

To show the benefit of including turning motion constraints, crossing speed, fuel consumption and traveling time are compared, as shown in Fig. 14 and Fig. 15. Fuel consumption and traveling time are defined as fuel and time consumed from 300 m before the intersection to 300 m after the intersection and reaching the original speed. If turning motion is not considered, the resultant fuel consumption is lower and the traveling time is less. On average, the fuel consumption is lower by 0.77% and the traveling time is lower by 6.00%. However, the crossing speed is higher the comfortable crossing speed. With the low speed limit constraint in this case, all the crossing speeds satisfy the safety constraint, however, none of them satisfy the riding comfort constraint. The crossing speed is higher by 39.32% on average for different traffic signal phase. If the speed limit is higher than the safety crossing speed, the safety crossing of the intersection cannot be assured.

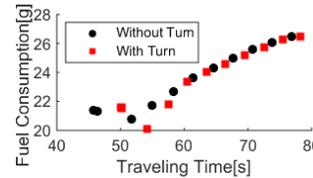
Fig. 14 Fuel Consumption Comparison for Turning Motion

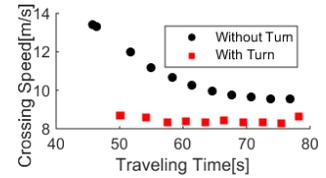
Fig. 15 Crossing Speed Comparison for Turning Motion

*C. Parametric Study of Weighting Parameters*

The simulation is carried out for both of single intersection and multiple intersections cases with a randomly generated traffic signal profile. The single intersections case is used to demonstrate the effect on fuel consumption and acceleration, and the multiple intersections case is used to demonstrate the influence on the intersection crossing window. The acceleration trajectories are shown in Fig. 16; the fuel consumption and traveling time results are shown in Fig. 17.

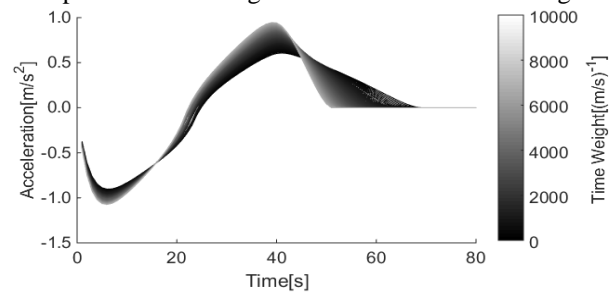
Fig. 16 Acceleration Trajectories for Different Time Weights

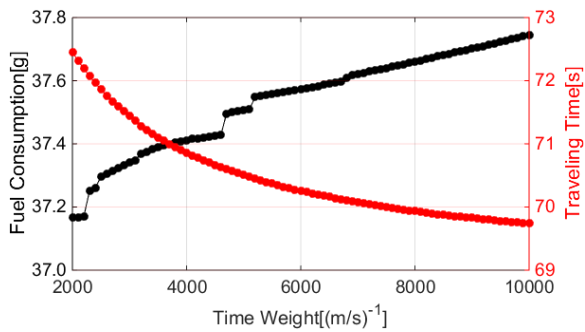

Fig. 17 Fuel Consumption and Traveling Time for Different Time Weights

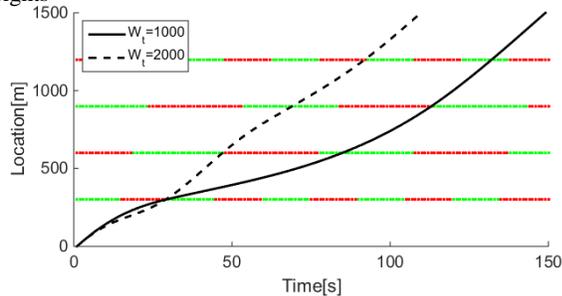

Fig. 18 Vehicle Trajectories for Different Time Weight

It can be seen from the motion trajectories that with an increasing weight of travel time, more aggressive acceleration is used to increase the average speed during the planning horizon. In addition, fuel consumption shows a non-linear increase with the increase in the time weighting parameter, while travel time shows polynomial decrease, due to the nonlinear nature of the fuel consumption function. The motion trajectories for multiple intersections case are shown in Fig. 18. It can be seen from the figures that with the increase in time weight, the vehicle tends to use more aggressive acceleration to catch earlier crossing windows, such that the fuel consumption during the planning horizon is increased.

## V. CONCLUSION

With the aid of broadcast traffic signal information, a vehicle's speed trajectory can be optimized for signalized intersections. We show that not only fuel consumption but also travel time can be reduced. However, the analysis is based on the ideal assumption that no other vehicle is present in that section of the road. For our next step, we will improve the robustness in eco-approach/departure by incorporating the proposed method into our previous studies on the maneuvers of other vehicles [23, 24]. In addition, the analysis is based on the connected automated vehicle assumption, which means the speed trajectory is followed precisely. However, a driver-assistance speed advisory would be a more practical application for the current development status of the automated vehicle. In this way, an analysis of drivers' responses to the optimal speed planning would also be a meaningful study.